\documentclass[a4paper,onecolumn]{revtex4}

\usepackage{amsmath}    
\usepackage{verbatim}   
\usepackage{bbold}
\usepackage{mathrsfs}
\usepackage{textcomp}
\usepackage{amsthm}
\usepackage{amssymb}
\usepackage{scrextend}
\usepackage{titlesec}
\usepackage{soul}

\usepackage{geometry}
	\usepackage{amsmath}
	\usepackage{amsfonts}
	\usepackage[ansinew]{inputenc}
	\usepackage[usenames,dvipsnames]{pstricks}
	\usepackage{hyperref}
	\usepackage{cleveref}
	\usepackage{scrextend}
	\usepackage{amsthm}
\changefontsizes[15pt]{11pt}

	\newtheorem{theorem}{Theorem}
	\newtheorem{lemma}{Lemma}
	\newtheorem*{theorem*}{Theorem}
	\newtheorem{corollary}{Corollary}

\allowdisplaybreaks

\begin{document}

\title{$q$-analogs of sinc sums and integrals}

\author{Martin Nicholson}

\begin{abstract}
$q$-analogs of sum equals integral relations $\sum_{n\in\mathbb{Z}}f(n)=\int_{-\infty}^\infty f(x)dx$ for sinc functions and binomial coefficients are studied. Such analogs are already known in the context of $q$-hypergeometric series. This paper deals with multibasic `fractional' generalizations that are not $q$-hypergeometric functions.
\end{abstract}

\maketitle

\section*{Introduction}
Surprizing properties of sinc sums and integals were first discovered by C. Stormer in 1895 [\onlinecite{stormer1},\onlinecite{stormer2}]. The more general properties of band limited functions were known to engineers from signal processing and to physicists. For example, K.S. Krishnan viewed them as a rich source for finding identities [\onlinecite{krishnan}]. R.P. Boas has studied the error term when approximating a sum of a band limited function with corresponding integral [\onlinecite{boas}]. More recently these properties were studied and popularized in a series of papers [\onlinecite{borwein1,borwein2,almkvist}].

$\mathrm{sinc}$ function is a special case of binomial coefficients
$$
\binom{2}{1+x}=\frac{\Gamma(3)}{\Gamma(1+x)\Gamma(1-x)}=\frac{2\sin\pi x}{\pi x}=2\,\mathrm{sinc}(\pi x).
$$
Therefore only sums with binomial coefficients will be studied in the following.
It is known that binomial coefficients are band limited (e.g., see [\onlinecite{pollard}])
$$
\binom{a}{u}=\frac{1}{2\pi}\int_{-\pi}^\pi (1+e^{it})^a e^{-iut}\,dt,
$$
i.e. their Fourier spectrum is limited to the band $|t|<\pi$. According to general theorems [\onlinecite{boas,borwein1}] whenever Fourier spectrum of a function $f(x)$ is limited to the band $|t|<2\pi$ one expects that
\begin{equation}\label{sumeqint}
\sum_{n=-\infty}^\infty f(n)=\int_{-\infty}^\infty f(x)dx.
\end{equation}
Bandwidth of a product of bandlimited functions is the sum of their bandwidths [\onlinecite{almkvist}]. In case of binomial coefficients this together with the theorem mentioned above implies that
\begin{equation}\label{binomsum}
	\sum_{n=-\infty}^\infty \binom{a}{\alpha n}^l=\int_{-\infty}^\infty \binom{a}{\alpha x}^ldx,\qquad 0<\alpha\le \frac{2}{l}.
\end{equation}
For a general band limited function the above formula would have been valid only when $\alpha< \frac{2}{l}$. The validity of \eqref{binomsum} when $\alpha= \frac{2}{l}$ is explained by the fact that spectral density of binomial coefficient vanishes at boundary values $t=\pm \pi$.

$q$-analog of the Gamma function is defined as
$$
\Gamma_q(x)=\frac{(q;q)_\infty}{(q^x;q)_\infty}(1-q)^{1-x}
$$
and the $q$-binomial coefficients
$$
\left[{a\atop b}\right]_q=\frac{\Gamma_q(a+1)}{\Gamma_q(b+1)\Gamma_q(a-b+1)},
$$
with the standard notations for the $q$-shifted factorials
$$
(a;q)_n=\prod_{k=0}^{n-1} (1-aq^k),\qquad (a_1,\ldots,a_r;q)_n=\prod_{k=1}^r(a_k;q)_n,\qquad (a;q)_\infty=\prod_{k=0}^{\infty} (1-aq^k).
$$
In the limit $q\to 1^{-}$ one has $\Gamma_q(a)\to \Gamma(a)$, i.e. standard values of the Gamma function and binomial coefficients are recovered.[\onlinecite{gasper}]

$q$-analog of the property of bandlimitedness has been studied in the literature [\onlinecite{ismail}]. This paper has a much more narrow scope and only deals with sums of binomial coefficients. We will find that \eqref{binomsum} with $0<\alpha\le 1/l$ has a very natural $q$-analog given in Theorem \ref{main}. However no such simple direct $q$-analog of \eqref{binomsum} with $1/l<\alpha\le 2/l$ is known. Nevertheless there is a formula that in the limit $q\to 1^{-}$ can be brought to the form \eqref{binomsum} after a series of simple steps which will be given in Theorem \ref{general}.

\section*{Main formula and its proof}

In the following, we will use a method of functional equations [\onlinecite{andrews}] (see also [\onlinecite{gasper}], sec. 5.2) combined with an idea due to G. Gasper [\onlinecite{gasper3}] to find a Laurent series expansion for a certain integral of an infinite product. First we need the following theorem taken from the book [\onlinecite{sidorov}]:
\begin{theorem}
Let 
\begin{equation}\label{F}
F(z)=\int_\gamma f(\zeta,z)d\zeta,
\end{equation}
where the following conditions are satisfied

{\rm{(1)}} $\gamma$ is an infinite picewise continous curve

{\rm{(2)}} the function $f(\zeta,z)$ is continous in $(\zeta,z)$ at $\zeta\in\gamma,~z\in D$, where $D$ is a domain in the complex $z$ plane,

{\rm{(3)}} for each fixed $\zeta\in\gamma$ the function $f(\zeta,z)$ viewed as a function of $z$ is regular in $D$,

{\rm{(4)}} integral \eqref{F} converges uniformly in $z\in D'$, where $D'$ is an arbitrary closed subdomain of $D$.

Then $F(z)$ is regular in $D$.
\end{theorem}

\begin{lemma}\label{L1}
Let  $p$ and $q$ be two real numbers that satisfy $0<p<q<1$, then
$$F(z)=\int_{-\infty}^\infty\frac{\left(bq^{\zeta},aq^{-\zeta};p\right)_\infty}{\left(-zq^{\zeta},-q^{1-\zeta}/z;q\right)_\infty}\,d\zeta$$
is regular in the half plane ${\mathrm{Re}}\, z>0$.
\end{lemma}
\noindent
$\it{Proof.}$ 
Put in the theorem above $f(\zeta,z)=\frac{\left(bq^{\zeta},aq^{-\zeta};p\right)_\infty}{\left(-zq^{\zeta},-q^{1-\zeta}/z;q\right)_\infty}$, $\gamma=(-\infty,+\infty)$, and $D$ an arbitrary domain in the half plane $\text{Re}~z > 0$. Then (1),(2) and (3) are obviously satisfied. To prove (4) let $p=e^{-\omega}$, $~q=p^\alpha$,$~\omega>0$,$~0 < \alpha < 1$ and consider the asymptotics of $f(\zeta,z)$ when $\zeta\to+\infty$. In this limit one has $(bq^\zeta;p)_\infty\to 1,~(-zq^\zeta;q)_\infty\to1$. According to an asymptotic formula ([\onlinecite{gasper}], p. 118)
$$
\mathrm{Re}[\ln(p^s;p)_\infty]=\frac{\omega}{2}(\mathrm{Re}\, s)^2+\frac{\omega}{2}(\mathrm{Re}\, s)+O(1),\qquad \mathrm{Re}\, s\to-\infty,
$$
we have
$$
|(aq^{-\zeta};p)_\infty|=|(p^{-\alpha\zeta-\omega^{-1}\ln a};p)_\infty|=O\left(|a|^{\alpha\zeta} q^{-(\alpha \zeta^2-\zeta)/2}\right),
$$
$$
|(-q^{1-\zeta}/z;q)_\infty|=|(q^{1-\zeta+\alpha^{-1}\omega^{-1}\ln z};q)_\infty|=O\left(|q/z|^{\zeta} q^{-(\zeta^2-\zeta)/2}\right).
$$
So
$$
f(\zeta,z)=O\left(|{z}{a^\alpha}/q|^{\zeta} q^{(1-\alpha)\zeta^2/2}\right), \quad \zeta\to+\infty.
$$
Similarly
$$
f(\zeta,z)=O\left(|{b^\alpha}/z|^{-\zeta} q^{(1-\alpha)\zeta^2/2}\right), \quad \zeta\to-\infty.
$$
It is now easy to see that the integral (*) converges. Hence according to Weierstrass M-Test integral $F(z)$ converges uniformly in $z$ when $\mathrm{Re}\,z\ge\delta>0$.
As a result the function
$$
f(a,b,z)=\frac{(-z,-q/z;q)_\infty}{\ln\frac{1}{q}}\int\limits_0^\infty\frac{\left(bt/z,pz/at;p\right)_\infty}{\left(-t,-q/t;q\right)_\infty}\frac{dt}{t}
$$
is regular when $\mathrm{Re}\,z > 0$\qed

\begin{lemma}\label{L2}
The function
$$
f(a,b,z)=\frac{(-z,-q/z;q)_\infty}{\ln\frac{1}{q}}\int_0^\infty\frac{\left(bt,a/t;p\right)_\infty}{\left(-zt,-q/(zt);q\right)_\infty}\frac{dt}{t}
$$
satisfies the functional equations
\begin{equation}\label{functional1}
f(a,b,z)=f(a,bp,z)-bf(a,bp,qz),
\end{equation}
\begin{equation}\label{functional2}
f(a,b,z)=f(ap,b,z)-af(ap,b,z/q).
\end{equation}
\end{lemma}
\noindent
$\it{Proof.}$ After a series of simple manipulations of the infinite products we find
\begin{align*}
f(a,b,qz)&=\frac{(-qz,-1/z;q)_\infty}{\ln\frac{1}{q}}\int_0^\infty\frac{\left(bt,a/t;p\right)_\infty}{\left(-qzt,-1/(zt);q\right)_\infty}\,\frac{dt}{t}\\
&=\frac{(-z,-q/z;q)_\infty}{z\ln\frac{1}{q}}\int_0^\infty \frac{z\left(bt,a/t;p\right)_\infty}{\left(-zt,-q/(zt);q\right)_\infty}\,dt\\
&=\frac{p(-z,-q/z;q)_\infty}{b\ln\frac{1}{q}}\int_0^\infty\frac{bt}{p}\frac{\left(bt,a/t;p\right)_\infty}{\left(-zt,-q/(zt);q\right)_\infty}\frac{dt}{t}\\
&=\frac{p}{b}(f(a,b,z)-f(a,b/p,z)).
\end{align*}
This is equivalent to \eqref{functional1}.
Similarly or using the first functional equation and the formula $f(a,b,z)=f(b,a,q/z)$ we find
\begin{align*}
f(a,b,z)=f(b,a,q/z)&=f(b,ap,q/z)-af(b,ap,q^2/z)\\
&=f(ap,b,z)-af(ap,b,z/q),
\end{align*}
as required. \qed

\begin{theorem}\label{main}
Let  $p$ and $q$ be two complex numbers such that $|p|<|q|<1$, then
\end{theorem}
$$\sum_{n=-\infty}^\infty (bq^n,aq^{-n};p)_\infty
    z^n q^{n(n-1)/2}=\frac{(-z,-q/z;q)_\infty}{\ln\frac{1}{q}}\int_0^\infty\frac{\left(bt/z,az/t;p\right)_\infty}{\left(-t,-q/t;q\right)_\infty}\frac{dt}{t}.
$$
$\it{Proof.}$ First consider the case $0 < p < q < 1$.  
The function $f(a,b,z)$ from Lemma \ref{L2} can be written in the form
$$
f(a,b,z)=(-z,-q/z;q)_\infty\int_{-\infty}^\infty\frac{\left(bq^{\zeta}/z,azq^{-\zeta};p\right)_\infty}{\left(-q^{\zeta},-q^{1-\zeta};q\right)_\infty}\,d\zeta.
$$
According to Lemma \ref{L1} $f(a,b,z)$ is a regular function of $z$ in the region $\mathrm{Re} z>0$. As a result $f(a,b,z)$ has the Laurent series expansion
$$
f(a,b,z)=\sum_{n=-\infty}^\infty c_n(a,b)z^n,\qquad \mathrm{Re}\, z>0.
$$
Functional equation \eqref{functional1} gives the following recursion relation for coefficients $c_n(a,b)$
$$
c_n(a,b)=(1-bq^n)c_n(a,bp).
$$
This recursion means that 
$$
c_n(a,b)=(bq^n;p)_\infty c_n(a,0).
$$
The functional equation \eqref{functional2} gives
$$
c_n(a,b)=(1-aq^{-n})c_n(a/p,b),
$$
from which one obtains 
$$
c_n(a,b)=(aq^{-n};p)_\infty c_n(0,b).
$$
By combining these equations one gets
$$
c_n(a,b)=(bq^n;p)_\infty c_n(a,0)=(bq^n,aq^{-n};p)_\infty c_n(0,0).
$$
It is known that ([\onlinecite{gasper}], ex. 6.16)
$$
\int_0^\infty\frac{1}{\left(-t,-q/t;q\right)_\infty}\frac{dt}{t}=(q;q)_\infty\ln\frac{1}{q}.
$$
According to Jacobi triple product formula
$$
(q,-z,-q/z;q)_\infty=\sum_{n=-\infty}^\infty
    z^n q^{n(n-1)/2}
$$
this implies that $c_n(0,0)=z^nq^{n(n-1)/2}$, so finally
$$
c_n(a,b)=(bq^n,aq^{-n};p)_\infty z^nq^{n(n-1)/2}.
$$
Now one needs to continue the result established for $\text{Re}~z > 0, 0 < p < q <1$ analytically to complex values of parameters $z,p,q$ to complete the proof. \qed

Series containing infinite products $(bq^n,aq^{-n};p)_\infty$ have been studied in [\onlinecite{ismail}]. It appears that the series in Theorem \ref{main} have been first considered in the paper [\onlinecite{vil}] which also contains a different representation for this sum in terms of an integral over a unit circle.

\section*{Consequences of the main formula}

\begin{corollary}\label{symmetric} The formula in Theorem \ref{main} can be written in symmetric form
$$\sum_{n=-\infty}^\infty \frac{(bq^n,aq^{-n};p)_\infty}{(-zq^n,-q^{1-n}/z;q)_\infty}
   =\int_{-\infty}^\infty\frac{\left(bq^x,aq^{-x};p\right)_\infty}{\left(-zq^x,-q^{1-x}/z;q\right)_\infty}\,dx,
$$
or in terms of $q$-binomial coefficients
\begin{equation}\label{binom}
\sum_{n=-\infty}^\infty \left[{a\atop b+\alpha n}\right]_{p}\frac{1}{(-zq^n,-q^{1-n}/z;q)_\infty}
   =\int_{-\infty}^\infty\left[{a\atop b+\alpha x}\right]_{p}\frac{1}{\left(-zq^x,-q^{1-x}/z;q\right)_\infty}\,dx,
\end{equation}
where $q=p^\alpha$, $0<\alpha<1$.
\end{corollary}

This gives an example of function for which sum equals integral. The case $|p|=|q| < 1$, $|b/a| < |z| < 1$ was known to Ramanujan. In this case, the series is Ramanujan's ${}_1\psi_1$ sum and the integral is Ramanujan's $q$-beta integral ([\onlinecite{gasper}], chs. 5,6).

Now let $z=e^{i\theta}$, $|\theta| < \pi$. Then
$$
\lim_{~\,q\to 1^{-}}\frac{\left(-z,-q/z;q\right)_\infty}{\left(-zq^x,-q^{1-x}/z;q\right)_\infty}=(1+z)^x(1+1/z)^{-x}=z^x.
$$
Let $q\to 1^{-}$ with $0<\alpha<1$ fixed in equation \eqref{binom}. Then formally
\begin{equation}\label{osler1}
\sum_{n=-\infty}^\infty \left({a\atop b+\alpha n}\right)e^{i\theta n}
   =\int_{-\infty}^\infty\left({a\atop b+\alpha x}\right) e^{i\theta x}\,dx,\quad ~0<\alpha<1.
\end{equation}
The range of validity of \eqref{osler1} is $-\pi\alpha<\theta < \pi\alpha$ as in \eqref{osler}, and not $-\pi<\theta < \pi$. Continuing formal manipulations we obtain by using \eqref{osler1} and binomial theorem
\begin{align}\label{osler2}
\nonumber\int_{-\infty}^\infty\left({a\atop b+\alpha x}\right) e^{i\theta x}\,dx&=\frac{1}{\alpha}e^{-i\theta b/\alpha}\int_{-\infty}^\infty\left({a\atop x}\right) e^{i\theta x/\alpha}\,dx\\
\nonumber&=\frac{1}{\alpha}e^{-i\theta b/\alpha}\sum_{n=-\infty}^\infty \left({a\atop n}\right)e^{i\theta n/\alpha}\\
\nonumber&=\frac{1}{\alpha}e^{-i\theta b/\alpha}\sum_{n=0}^\infty \left({a\atop n}\right)e^{i\theta n/\alpha}\\
&=\frac{1}{\alpha}e^{-i\theta b/\alpha}(1+e^{i\theta/\alpha})^a,\quad -\pi\alpha<\theta < \pi\alpha.
\end{align}
Finally \eqref{osler1} and \eqref{osler2} imply
\begin{equation}\label{osler}
\sum_{n=-\infty}^\infty \left({a\atop b+\alpha n}\right)v^{b+\alpha n}=\frac{1}{\alpha}(1+v)^a,\quad |v|=1,~ |\arg v|<\pi,~ 0<\alpha\le 1 ,
\end{equation}
which is T. Osler's generalization of binomial theorem [\onlinecite{osler}]. According to Osler [\onlinecite{osler}], the special case $\alpha=1$ of \eqref{osler} was first stated by Riemann [\onlinecite{riemann}]. It also follows from Ramanujan's ${}_1\psi_1$ sum in the limit $q\to 1^{-}$. 

It should be noted that while \eqref{osler} has a closed form, the series in Theorem \ref{main} does not. If $p=q^2$,$z=1$,$b=aq^2$, then one can prove that
$$
\sum_{n=-\infty}^\infty (bq^n,p/aq^n;p)_\infty
    z^n q^{n(n-1)/2}=2 \left(q a,q/a;q^2\right){}_{\infty }\sum _{n=-\infty}^\infty \frac{\left(-1/a\right)^n q^{n^2+n}}{1-a q^{2 n+1}}.
$$
The sum on the RHS is proportional to Appell-Lerch sum $m(qa^2,q^2,q^2/a)$ in the notation of the paper [\onlinecite{mortenson}]. In general Appell-Lerch sums do not have an infinite product representation. For example, by taking $a=q^{-1/2}$ in $m(qa^2,q^2,q^2/a)$ we get the sum of the type $m(1,q^2,z)$ which is related to mock theta function of order 2 (see formula (4.2) in [\onlinecite{mortenson}]). 

\begin{corollary}\label{main2}
The series
$$
\sum_{n=-\infty}^\infty \frac{(bq^n,p/aq^n;p)_\infty}{(-zq^n,-q/zq^n;q)_\infty},\qquad |p|<|q|
$$
with $p$ and $q$ fixed depends only on $b/z$ and $az$.
\end{corollary}

\begin{theorem}\label{fourier}
\begin{align*}
&\int_{-\infty}^\infty\frac{\left(bq^x,aq^{-x};p\right)_\infty}{\left(-q^x,-q^{1-x};q\right)_\infty}\,e^{ixy}\,dx\\
&=\frac{2\pi i/\log q}{\sinh\frac{\pi y}{\log q}}\frac{(-q,-q,e^{iy},qe^{-iy};q)_\infty}{(q,q,-e^{iy},-qe^{-iy};q)_\infty}\sum_{n=-\infty}^\infty \frac{(bq^n,aq^{-n};p)_\infty}{(-q^n,-q^{1-n};q)_\infty}\, e^{iny}.
\end{align*}
\end{theorem}
\noindent
$\it{Proof.}$ Consider the contour integral
$$
\int_C\frac{\left(bq^z,aq^{-z};p\right)_\infty}{\left(-q^z,-q^{1-z};q\right)_\infty}\,e^{izy}\,dz
$$
where $C$ is rectangle with vertices at $(\pm R,0),~(\pm R,-{2\pi i}/{\log q})$. In view of asymptotics found in the proof of Lemma \ref{L1} integrals over the vertical segments vanish in the limit $R\to +\infty$. Integrals over the horizontal segments are convergent and related by a factor of $-e^{2\pi y/\log q}$. The integrand has simple poles at $z=n-{\pi i}/{\log q}$ with residues
$$
-\frac{e^{\pi y/\log q}}{(q;q)_\infty^2\log q}(-b q^n,-a q^{-n};p){}_{\infty }(-1)^n q^{n(n-1)/2} e^{i n y}.
$$
Application of the residue theorem yields
$$ 
\int_{-\infty}^\infty\frac{\left(bq^x,aq^{-x};p\right)_\infty}{\left(-q^x,-q^{1-x};q\right)_\infty}\,e^{ixy}\,dx=\frac{\pi  i/\log q}{ (q;q)^2_{\infty }\sinh \frac{\pi  y}{\log q}}\sum _{n=-\infty}^\infty (-b q^n,-a q^{-n};p){}_{\infty }(-1)^n q^{n(n-1)/2} e^{i n y}.
$$
According to Corollary \ref{main2}
$$
\sum _{n=-\infty}^\infty (-b q^n,-a q^{-n};p){}_{\infty }(-1)^n q^{n(n-1)/2} e^{i n y}=\frac{(e^{iy},qe^{-iy};q)_\infty}{(-e^{iy},-qe^{-iy};q)_\infty}\sum _{n=-\infty}^\infty (b q^n,a q^{-n};p){}_{\infty } q^{n(n-1)/2} e^{i n y}.
$$
To complete the proof observe that
$$
\sum _{n=-\infty}^\infty (b q^n,a q^{-n};p){}_{\infty } q^{n(n-1)/2} e^{i n y}=(-1,-q;q)_\infty\sum_{n=-\infty}^\infty \frac{(bq^n,aq^{-n};p)_\infty}{(-q^n,-q^{1-n};q)_\infty}\, e^{iny}
$$
and $(-1,-q;q)_\infty=2(-q;q)^2_\infty$.\qed

One can see from Theorem \ref{fourier} that the function 
$$
g(x)=\frac{\left(bq^x,aq^{-x};p\right)_\infty}{\left(-q^x,-q^{1-x};q\right)_\infty}
$$
is not band limited. However Fourier transform of $g(x)$ vanishes at frequencies $y=2\pi m$, where $m\neq 0$ is an integer. Hence according to Poisson summation formula [\onlinecite{titchmarsh}]
$$
\sum _{n=-\infty}^\infty g(x)=\sum _{n=-\infty}^\infty\int_{-\infty}^\infty g(x)e^{-2\pi in x}dx=\int_{-\infty}^\infty g(x)dx
$$
in agreement with Corollary \ref{symmetric}.

The fact that bilateral summation formulas in the theory of $q$-hypergeometric functions give examples of functions of the type \eqref{sumeqint} has been recognized in the literature. 

\begin{corollary}
Let $|p|<|q|$ and $m\in\mathbb{Z}$, then
$$
\int_{-\infty}^\infty\frac{\left(bq^x,aq^{-x};p\right)_\infty}{\left(-q^x,-q^{1-x};q\right)_\infty}\,q^{mx}\,dx=\sum_{n=-\infty}^\infty \frac{(bq^n,aq^{-n};p)_\infty}{(-q^n,-q^{1-n};q)_\infty}\, q^{mn}.
$$
\end{corollary}
\noindent
$\it{Proof.}$ Resolve the $\frac{0}{0}$ ambiguity at the rhs of the formula of Theorem 2 using L'Hopital's Rule.\qed

Next we apply the method due to Bailey [\onlinecite{bailey}] to the identity in Theorem 2.
\begin{theorem}\label{bailey}
$$
\sum _{n=-\infty}^\infty  \left(b_1q^n,b_2q^n,a_1 q^{-n},a_2q^{-n};p\right)_{\infty }z^nq^{n (n-1)}=z\sum _{n=-\infty}^\infty  \left(b_1q^n/z,b_2q^n/z,a_1zq^{-n},a_2zq^{-n};p\right)_{\infty }z^{-n}q^{n (n-1)}.
$$
\end{theorem}
\noindent
$\it{Proof.}$ Multiplying the equations
$$
\sum_{n=-\infty}^\infty (b_1q^n,a_1q^{-n};p)_\infty
    e^{i\theta n} q^{n(n-1)/2}=\frac{(-e^{i\theta},-qe^{-i\theta};q)_\infty}{\ln\frac{1}{q}}\int_0^\infty\frac{\left(b_1te^{-i\theta},a_1e^{i\theta}/t;p\right)_\infty}{\left(-t,-q/t;q\right)_\infty}\frac{dt}{t},
$$
$$
\sum_{n=-\infty}^\infty (b_2q^n,a_2q^{-n};p)_\infty
    e^{-i\theta n}z^n q^{n(n-1)/2}=\frac{(-ze^{-i\theta},-qe^{i\theta}/z;q)_\infty}{\ln\frac{1}{q}}\int_0^\infty\frac{\left(b_2te^{i\theta}/z,a_2ze^{-i\theta}/t;p\right)_\infty}{\left(-t,-q/t;q\right)_\infty}\frac{dt}{t},
$$
and integrating with respect to $\theta$ one obtains
\begin{align*}
&\sum _{n=-\infty}^\infty  \left(b_1q^n,b_2q^n,a_1 q^{-n},a_2q^{-n};p\right)_{\infty }z^nq^{n (n-1)}\\
&=\int_{-\pi}^\pi\frac{d\theta}{2\pi}\frac{(-e^{i\theta},-qe^{-i\theta};q)_\infty}{\ln\frac{1}{q}}\int_0^\infty\frac{\left(b_1t_1e^{-i\theta},a_1e^{i\theta}/t_1;p\right)_\infty}{\left(-t_1,-q/t_1;q\right)_\infty}\frac{dt_1}{t_1}\\
&\times\frac{(-ze^{-i\theta},-qe^{i\theta}/z;q)_\infty}{\ln\frac{1}{q}}\int_0^\infty\frac{\left(b_2t_2e^{i\theta}/z,a_2ze^{-i\theta}/t_2;p\right)_\infty}{\left(-t_2,-q/t_2;q\right)_\infty}\frac{dt_2}{t_2} \\
&=z\int_{-\pi}^\pi\frac{d\theta}{2\pi}\frac{(-e^{-i\theta},-qe^{i\theta};q)_\infty}{\ln\frac{1}{q}}\int_0^\infty\frac{\left(b_2t_2e^{i\theta}/z,a_2ze^{-i\theta}/t_2;p\right)_\infty}{\left(-t_2,-q/t_2;q\right)_\infty}\frac{dt_2}{t_2}\\
&\times\frac{(-e^{i\theta}/z,-qze^{-i\theta};q)_\infty}{\ln\frac{1}{q}}\int_0^\infty\frac{\left(b_1t_1e^{-i\theta},a_1e^{i\theta}/t_1;p\right)_\infty}{\left(-t_1,-q/t_1;q\right)_\infty}\frac{dt_1}{t_1}\\
&=z\sum _{n=-\infty}^\infty  \left(b_1q^n/z,b_2q^n/z,a_1zq^{-n},a_2zq^{-n};p\right)_{\infty }z^{-n}q^{n (n-1)}.\qed
\end{align*}

It is straightforward to rewrite the formula of Theorem \ref{bailey} in terms of $q$-binomial coefficients:
\begin{corollary} Let $0<q<1$ and $0<\alpha<1$, then
$$
\sum_{n=-\infty}^\infty\left[{a_1\atop b_1+\alpha n}\right]_p\left[{a_2\atop b_2+\alpha n}\right]_pp^{\alpha n(n-1)+\theta n}=p^\theta\sum_{n=-\infty}^\infty\left[{a_1\atop b_1-\theta+\alpha n}\right]_p\left[{a_2\atop b_2-\theta+\alpha n}\right]_pp^{\alpha n(n-1)-\theta n}.
$$
\end{corollary}

Theorem \ref{main} can be generalized.
\begin{theorem}\label{general}
Let  $q=p_1^{\alpha_1}=p_2^{\alpha_2}$ where $\alpha_1>0$, $\alpha_2>0$, and $\alpha_1+\alpha_2<1$. Then
\end{theorem}
\begin{align*}
\sum_{n=-\infty}^\infty &\left[{a_1\atop b_1+\alpha n}\right]_{p_1}\left[{a_2\atop b_2+\alpha n}\right]_{p_2}\frac{1}{(-zq^n,-q^{1-n}/z;q)_\infty}\\&=\int_{-\infty}^\infty\left[{a_1\atop b_1+\alpha x}\right]_{p_1}\left[{a_2\atop b_2+\alpha x}\right]_{p_2}\frac{dx}{\left(-zq^x,-q^{1-x}/z;q\right)_\infty}.
\end{align*}

In fact, the number of $q$-binomial coefficients in this formula can be arbitrary as long as the parameters $\alpha_j>0$ are subject to the constraint $\sum_j\alpha_j<1$.


\begin{thebibliography}{9}

\bibitem{stormer1} C. St\"ormer, {\it{Om en generalisation af integralet}} $\int_0^\infty\frac{\sin{ax}}{x}dx=\frac{\pi}{2}$, Videnskapsselskapet i Kristiania, {\bf{4}}, (1895).

\bibitem{stormer2} C. St\"ormer, {\it{Sur une G\'en\'eralisation de la formule}} $\frac{\varphi }{2} = \frac{{\sin  \varphi }}{1} - \frac{{\sin  2\varphi }}{2} + \frac{{\sin  3\varphi }}{3}\ldots $, Acta math. {\bf{19}}, 341-350 (1895).

\bibitem{krishnan} K S Krishnan, {\it{On the equivalence of certain infinite series and the corresponding integrals}}, J. Indian Math. Soc., {\bf{12}}, 79-88, (1948).

\bibitem{simon} R. Simon, {\it{K.S. Krishnan's 1948 Perception of the Sampling Theorem}}, Resonance, {\bf{7}}, 20-25 (2002).

\bibitem{boas} R.P. Boas, {\it{Summation formulas and band-limited signals}}, Tohoku Math. Journ. {\bf{24}}, 121-125 (1972).

\bibitem{borwein1} D. Borwein and J.M. Borwein, {\it{Some remarkable properties of sinc and related integrals,}}
The Ramanujan Journal, {\bf{5}}, 73-89 (2001).

\bibitem{borwein2} R. Baillie, D. Borwein and J.M. Borwein, {\it{Surprising sinc sums and integrals}}, Amer.
Math. Monthly {\bf{115}}, 888-901 (2008).

\bibitem{almkvist} G. Almkvist, J. Gustavsson, {\it{More remarkable sinc integrals and sums}}, \href{https://arxiv.org/abs/1405.1265}{arXiv:1405.1265} (2014).

\bibitem{ismail} M.E. Ismail, A.I. Zayed, $q${\it{-analogue of the Whittaker-Shannon-Kotel'nikov sampling theorem}}, Proc. Amer. Math. Soc. {\bf{131}}, 3711-3719 (2003).

\bibitem{pollard} H. Pollard, O. Shisha, {\it{Variations on the Binomial Series }}, Amer.
Math. Monthly {\bf{79}}, 495-499 (1972).

\bibitem{gasper} G. Gasper, M. Rahman, {\it{Basic hypergeometric series}}, 2nd ed., Cambridge University Press, Cambridge (2004).

\bibitem{ismail} M.E.H. Ismail, D. Stanton, {\it{$q$-Taylor theorems, polynomial expansions, and
interpolation of entire functions}}, J. Approx. Theory {\bf{123}}, 125-146 (2003).

\bibitem{andrews} G.E. Andrews, R. Askey, {\it{A simple proof of Ramanujan's summation of the
${}_1\psi_1$}}, Aequationes Math. {\bf{18}}, 333-337 (1978).

\bibitem{gasper3} G. Gasper, Solution to problem 6497 ($q$-Analogues of a gamma function identity,
by R. Askey), Amer. Math. Monthly {\bf{94}}, 199-201 (1987).

\bibitem{sidorov} Yu.V. Sidorov, M.V. Fedoryuk, and M.I. Shabunin, {\it{Lectures on the Theory of Functions of a Complex Variable}} (English translation), Mir Publishers, Moscow (1985).

\bibitem{rahman} M.E.H. Ismail and M. Rahman, {\it{Some bilateral sums and integrals}}, Pacific J. Math. {\bf{170}}, 497-515 (1995).

\bibitem{vil} N.M. Vildanov, {\it{Some extensions of Ramanujan's ${}_1\psi_1$ summation formula}}, \href{https://arxiv.org/abs/1204.6569}{arXiv:1204.6569 } (2012).

\bibitem{osler} T.J. Osler, {\it{Taylor's series generalized for fractional derivatives and applications}}, SIAM J. Math. Anal. {\bf{2}}, 37-48 (1971).

\bibitem{mortenson} D.R. Hickerson, E.T. Mortenson, {\it{Hecke-type double sums, Appell-Lerch sums, and mock theta functions, I}}, Proc. Lond. Math. Soc. {\bf{109}}, 382-422 (2014).

\bibitem{titchmarsh} E.C. Titchmarsh, {\it{Introduction to the Theory of Fourier Integrals}},  2nd.ed., Oxford University Press (1948).

\bibitem{suslov} S.K. Suslov, {\it{Multiparameter Ramanujan-Type $q$-Beta Integrals}}, Ramanujan J. {\bf{2}}, 351-369 (1998).

\bibitem{bailey} W.N. Bailey, {\it{On the basic bilateral hypergeometric series}} ${}_2\psi_2$, Quart. J. Math.
{\bf{1}}, 194-198 (1950).

\bibitem{butzer} P.L. Butzer, P.J.S.G. Ferreira, J.R. Higgins, S. Saitoh, G. Schmeisser, R.L. Stens, {\it{Interpolation and Sampling: E.T. Whittaker, K. Ogura and Their Followers}}, J. Fourier Anal. Appl. {\bf{17}}, 320-354 (2011).

\bibitem{riemann} B. Riemann, {\it{Versuch einer allgemeinen Auffasung der Integration und Differentiation}}, The Collected
Works of Bernhard Riemann, 2nd ed., H. Weber, ed., Dover, New York, (1953). pp. 353-366

\bibitem{gasper2} G. Gasper, {\it{Elementary derivations of summation and transformation formulas for
$q$-series}}, Fields Inst. Commun. {\bf{14}}, 55-70. (1997).

\end{thebibliography}
\end{document}